\documentclass[a4paper,leqno,11pt]{amsart}

\usepackage{amsmath}
\usepackage{amscd}
\usepackage{amssymb}
\usepackage{enumerate}
\usepackage{indentfirst}
\usepackage{latexsym}
\usepackage{multicol}
\usepackage{amsthm}

\newcommand{\Int}{\displaystyle \int}

\numberwithin{equation}{section}
\def\div{ \hbox{\rm div}\,  }

\def\Id{\hbox{\rm Id}}
\def\adj{\hbox{\rm adj}}
\def\Tr{\hbox{\rm Tr}}
\newcommand{\R}{{\mathbb R}}
\newcommand{\N}{{\mathbb N}}

\newcommand{\Z}{{\mathbb Z}}

\def\d{\partial}

\def\ddj{\dot\Delta_j}

\def\cM{{\mathcal M}}
\def\cC{{\mathcal C}}
\def\cS{{\mathcal S}}
\let\tilde=\widetilde

\newcommand{\dC}{\delta\! C}
\newcommand{\df}{\delta\! f}
\newcommand{\dg}{\delta\! g}
\newcommand{\du}{\delta\! u}
\newcommand{\dv}{\delta\! v}
\newcommand{\dU}{\delta\! U}
\newcommand{\dP}{\delta\!P}
\newcommand{\dQ}{\delta\!Q}
\newcommand{\dR}{\delta\!R}
\newcommand{\dr}{\delta\!\rho}
\newcommand\dPi{\delta\!\Pi}

\newtheorem{lem}{Lemma}
\newtheorem{cor}{Corollary}
\newtheorem{prop}{Proposition}
\newtheorem{theo}{Theorem}
\newtheorem{rem}{Remark}

\newenvironment{p}{
\noindent\textit{\textbf{Proof:}}~}
{\hfill\rule{2.1mm}{2.1mm}
}

\begin{document}
\title[Inhomogeneous incompressible Navier-Stokes equations]{A Lagrangian approach for the incompressible Navier-Stokes equations
with variable density}

\author[R. Danchin]{Rapha\"el Danchin}
\address[R. Danchin]
{Universit\'e Paris-Est, LAMA, UMR 8050,
 61 avenue du G\'en\'eral de Gaulle,
94010 Cr\'eteil Cedex, France.}
\email{danchin@univ-paris12.fr}

\author[P.B. Mucha]{Piotr Bogus\l aw Mucha}
\address[P.B. Mucha]{Instytut Matematyki Stosowanej i Mechaniki,
 Uniwersytet Warszawski, 
ul. Banacha 2,  02-097 Warszawa, Poland.} 
\email{p.mucha@mimuw.edu.pl}
\date\today

\begin{abstract} Here we  investigate the Cauchy problem for  the inhomogeneous Navier-Stokes equations in the whole $n$-dimensional space. Under
some smallness assumption on the data, we show
the existence of  global-in-time   unique solutions in a critical functional framework. The initial density is required to belong  to the multiplier space of 
$\dot B^{n/p-1}_{p,1}(\R^n)$. In particular, piecewise constant initial densities
are admissible data \emph{provided the jump at the interface is small enough}, and generate global unique solutions
with piecewise constant densities.  
Using Lagrangian coordinates is the key to our results as it enables us to solve the
system by means of the basic contraction mapping theorem. As a consequence, 
conditions for uniqueness are the same as for existence.
\end{abstract}
\maketitle



\section*{Introduction}

 We address the well-posedness issue for the the incompressible Navier-Stokes equations 
with variable density in the whole space $\R^n$: 
\begin{equation}\label{eq:euler}
\left\{\begin{array}{l}
\d_t\rho+u\cdot\nabla\rho=0,\\[1ex]
\rho(\d_tu+u\cdot\nabla u)-\mu\Delta u+\nabla P=0,\\[1ex]
\div u=0,\\[1ex]
u|_{|t=0}=u_0.
\end{array}\right.
\end{equation}
Above $\rho=\rho(t,x)\in\R_+$ stands for the density, 
$u=u(t,x)\in\R^N,$ for the velocity field and $P=P(t,x)\in\R,$ for 
the pressure field. The viscosity coefficient $\mu$ is a given positive
real number.
We supplement this system with the following boundary conditions: 
\begin{itemize}
\item the velocity $u$  tends to $0$ at infinity,
\item the density tends to some positive constant
$\rho^*$ at infinity.
\end{itemize} 
The exact meaning of those boundary conditions will be given by the functional framework
in which we shall solve the system. 
In what follows, we will  take $\rho^*=1$ to simplify the presentation.
\smallbreak

This old system has known a renewed interest recently, among the mathematics community. 
The existence of strong smooth solutions with positive density 
has been established in e.g. \cite{LS} 
whereas the theory of global weak solutions with finite energy has been 
performed in the book \cite{Lions} by P.-L. Lions (see also the references
therein, and the monograph \cite{AKM}). 
As pointed out in \cite{CK},
 it is possible to construct strong unique solutions for some classes of smooth enough data 
 with vanishing density. 
 
 In the present paper, we aim at solving the above system 
 in \emph{critical} functional spaces, that is  in spaces which have 
 the same invariance with respect to time and space dilation as 
 the system itself (see e.g. \cite{D1} or \cite{DM2} for more explanations
 about this nowadays classical approach). In this framework, it has been 
 stated in 
 \cite{Abidi,D1} that,  for data $(\rho_0,u_0)$ such that 
$$(\rho_0-1)\in\dot B^{n/p}_{p,1}(\R^n),\qquad
u_0\in\dot B^{n/p-1}_{p,1}(\R^n)\quad\hbox{with}\quad\div u_0=0$$ 
and that, for a small enough constant $c,$
\begin{equation}\label{eq:smalldata0}
\|\rho_0-1\|_{\dot B^{n/p}_{p,1}(\R^n)}+\mu^{-1}\|u_0\|_{\dot B^{n/p-1}_{p,1}}\leq c,
\end{equation}
we have for any  $p\in[1,2n)$: 
\begin{itemize}
\item existence of a global solution $(\rho,u,\nabla P)$ 
with  $\rho-1\in\cC_b(\R_+;\dot B^{n/p}_{p,1}),$
$u\!\in\!\cC_b(\R_+;\dot B^{n/p-1}_{p,1})$ and $\d_tu,\nabla^2u,\nabla P\in L_1(\R_+;\dot B^{n/p-1}_{p,1})$;
\item uniqueness in the above space if in addition  $p\leq n.$
\end{itemize}
Those results have been somewhat extended in \cite{AP} where it has been noticed 
that $\rho_0-1$ may be taken in a larger Besov space, with another Lebesgue exponent.
\smallbreak
The above  results are based on maximal regularity estimates in Besov spaces for the evolutionary Stokes system, 
and on the Schauder-Tychonoff fixed point theorem.
In effect, owing to the hyperbolicity of the density equation, 
there is a loss of one derivative in the stability estimates
thus precluding the use of the contraction mapping (or Banach fixed point) theorem. 
As a consequence, the conditions for uniqueness are \emph{stronger} than those
for  existence. 

Let  us also point out that all the above results concerning existence \emph{with} 
uniqueness require the density to be at least uniformly continuous. 
This condition has been somewhat weakened recently by P. Germain in \cite{Germain}. 
However, there, initial densities with jump across an interface cannot be considered.

\medbreak
In the present paper, we aim at solving System \eqref{eq:euler} in the Lagrangian coordinates. 
The main motivation is that the density is \emph{constant} along the flow so that
only the (parabolic type) equation for the velocity has to be considered. 
We shall show that, after performing this change of coordinates, solving
\eqref{eq:euler} may be done by means of the Banach fixed point theorem. 
As a consequence, the condition for uniqueness need not be stronger than
that for the existence, and  the flow map is Lipschitz continuous.

Our main result states the global in time existence of regular solutions to the inhomogeneous Navier-Stokes equations in $\R^n$ in
optimal Besov setting, under suitable smallness of the data. 
As regards the initial density, the admissible regularity is so low 
it may have (small) jumps
across a $C^1$  interface. 
 This is  of particular interest from the viewpoint of physics
 as it implies  that motion of a mixture of two incompressible fluids with slightly different
densities can be modeled by the inhomogeneous Navier-Stokes equations. In addition, the regularity of the 
constructed velocity  suffices to preserve  the 
$C^1$ regularity of the interface between fluids.

\medbreak
We now come to the plan  of the paper. 
In the next section, we present our main results and give some insight of the  proof. 
Section \ref{s:linear}  is devoted to solving the linearized system \eqref{eq:euler}
in Lagrangian coordinates. 
This will enable us to define a map $\Phi:E_p\rightarrow E_p$
where $E_p$ stands for  the functional space in which the Lagrangian version 
of the momentum equation of \eqref{eq:euler} is going to be solved. 
That  $\Phi$ fulfills the conditions of the contraction mapping
theorem on a small enough ball of $E_p$ is the main purpose of  Section \ref{s:nonlinear}.
  In the Appendix we prove several important  results concerning
the Lagrangian coordinates and Besov spaces.

\medbreak
\noindent {\bf Notation:}  
Throughout, the notation $C$ stands for a generic constant (the meaning of which depends on the context), and we sometimes write $X\lesssim Y$ instead of $X\leq CY.$ 
 Finally, for $A=(A_{ij})_{1\leq i,j\leq n}$ and  $B=(B_{ij})_{1\leq i,j\leq n}$
 two matrices, we denote $A:B=\sum_{i,j} A_{ij} B_{ji}.$


 \section{Main results and principle of the proof}

Let us first derive formally the Lagrangian equations
corresponding to \eqref{eq:euler}\footnote{The reader may refer to the Appendix for the rigorous derivation
in our functional setting.}.
Let $X_u$ be the flow associated to the vector-field $u,$ that is
the solution to 
\begin{equation}\label{lag}
X_u(t,y)=y+\int_0^tu(\tau,X_u(\tau,y))\,d\tau.
\end{equation}
Let us recall that by  Liouville's formula for transport equations,  the divergence-free condition 
is equivalent to  $|DX_u|\equiv1.$ In other words the map (\ref{lag})
is measure preserving.
Now, denoting 
$$
\bar\rho(t,y):=\rho(t,X_u(t,y)),\quad
\bar P(t,y):= P(t,X_u(t,y))\ \hbox{ and }\ 
\bar u(t,y)=u(t,X_u(t,y))
$$
with $(\rho,u,\nabla P)$ a solution of \eqref{eq:euler}, 
and using the chain rule and Lemma \ref{l:div} from the Appendix, we gather
that $\bar\rho(t,\cdot)\equiv\rho_0$ and that $(\bar u,\nabla \bar P)$ 
satisfies
 \begin{equation}\label{eq:lagrangian}
\left\{\begin{array}{l}
\rho_0\d_t\bar u-\mu\div_y (A_u{}^T\!A_u\nabla_y\bar u)+{}^T\!A_u\nabla_y\bar P=0\qquad\hbox{with}\quad A_u=(D_yX_u)^{-1},\\[1ex]
\div_y(A_u\bar u)=0.\end{array}\right.
\end{equation}
Motivated
by prior works (see e.g. \cite{Abidi,AP,D1,DM2}) 
we want to solve the above system in \emph{critical} homogeneous Besov spaces.
Let us recall that, for $1\leq p\leq\infty$ and $s\leq n/p,$ a tempered distribution $u$ 
over $\R^n$ 
belongs to the homogeneous Besov space $\dot B^s_{p,1}(\R^n)$ if 
$$
u=\sum_{j\in\Z}\ddj u\quad\hbox{in }\ \cS'(\R^n)
$$
and 
$$
\|u\|_{\dot B^s_{p,1}(\R^n)}:=\sum_{j\in\Z}2^{js}\|\ddj u\|_{L_p(\R^n)}<\infty.
$$
Here $(\ddj)_{j\in\Z}$ denotes a homogeneous dyadic resolution of unity 
in Fourier variables (see e.g. \cite{BCD}, Chap. 2 for more details). 
\smallbreak
Loosely speaking, a function belongs to $\dot B^s_{p,1}(\R^n)$
if it as $s$ derivatives in $L_p(\R^n).$ 
In the present paper, we shall make an extensive  use of the following classical properties:
\begin{itemize}
\item the Besov space $\dot B^{n/p}_{p,1}(\R^n)$ is a Banach algebra embedded in 
the set of continuous functions going to $0$ at infinity, whenever $1\leq p<\infty;$
\item the usual product maps $\dot B^{n/p-1}_{p,1}(\R^n)\times \dot B^{n/p}_{p,1}(\R^n)$
in $\dot B^{n/p-1}_{p,1}(\R^n)$ whenever $1\leq p<2n.$
\end{itemize}

{}From now on, we shall omit $\R^n$ in the notation for Besov spaces. 
We shall obtain the existence and uniqueness of a global solution $(\bar u,\nabla\bar P)$ 
for \eqref{eq:lagrangian} in the space
$$
E_p:=\bigl\{(\bar u,\nabla \bar P): \bar u\in\cC_b(\R_+;\dot B^{n/p-1}_{p,1}), 
\ \d_t\bar  u,\nabla ^2\bar u, \nabla\bar P\in L_1(\R_+;\dot B^{n/p-1}_{p,1})\bigr\},
$$
and we shall endow $E_p$ with the norm
$$
\|(\bar u,\nabla\bar P)\|_{E_p}:=\|\bar u\|_{L_\infty(\R_+;\dot B^{n/p-1}_{p,1})}+
\|\d_t\bar u,\mu\nabla^2\bar u,\nabla\bar P\|_{L_1(\R_+;\dot B^{n/p-1}_{p,1})}.
$$
 We shall also use the \emph{local} version $E_p(T)$ of $E_p,$
 pertaining to functions defined on $[0,T)\times\R^n.$ 
 Writing out the exact definition and the corresponding norm is
 left to the reader.
\smallbreak

The required regularity for the initial density $\rho_0$ is that it belongs to the \emph{multiplier
space}  $\cM(\dot B^{n/p-1}_{p,1})$ for $\dot B^{n/p-1}_{p,1},$
that is the set of those distributions $\rho_0$ such that $\psi\rho_0$
is in $\dot B^{n/p-1}_{p,1}$ whenever $\psi$ is in $\dot B^{n/p-1}_{p,1},$
endowed with the norm
\begin{equation}\label{eq:defmult}
\|\rho_0\|_{\cM(\dot B^{n/p-1}_{p,1})}:=\sup \|\psi\rho_0\|_{\dot B^{n/p-1}_{p,1}}
\end{equation}
where the supremum is taken over those functions $\psi$ in $\dot B^{n/p-1}_{p,1}$ with norm $1.$

Let us now state our main result.
\begin{theo}\label{th:lagrangian}
 Let $p\in[1,2n)$ and  $u_0$ be a divergence-free vector field in $\dot B^{n/p-1}_{p,1}(\R^n).$
 Assume that the initial density $\rho_0$ belongs to the multiplier space $\cM(\dot B^{n/p-1}_{p,1}).$
 There exists a constant $c$ depending only on $p$ and on $n$ such that if
 $$
 \|\rho_0-1\|_{\cM(\dot B^{n/p-1}_{p,1})}+\mu^{-1}\|u_0\|_{\dot B^{n/p-1}_{p,1}}\leq c
 $$
 then System \eqref{eq:lagrangian} has a unique global solution $(\bar u,\nabla\bar P)$ in $E_p.$
Moreover, we have
$$
\|(\bar u,\nabla\bar P)\leq C\|u_0\|_{\dot B^{n/p-1}_{p,1}}
$$ 
for some constant $C$ depending only on $n$ and on $p,$
and the flow map $(\rho_0,u_0)\longmapsto (\bar u,\nabla\bar P)$
is Lipschitz  continuous from  $\cM(\dot B^{n/p-1}_{p,1})\times\dot B^{n/p-1}_{p,1}$
to $E_p.$
 \end{theo}
 In the case where only the density satisfies the smallness condition, 
we get the following local-in-time existence result:
\begin{theo}\label{th:local}
Under the above regularity assumptions, there exists a constant $c$ depending only on $p$ and on $n$ such that if
 $$
 \|\rho_0-1\|_{\cM(\dot B^{n/p-1}_{p,1})}\leq c
 $$
 then there exists some $T>0$ such that 
 System \eqref{eq:lagrangian} has a unique local solution $(\bar u,\nabla\bar P)$ in $E_p(T),$
 and the flow map $(\rho_0,u_0)\longmapsto (\bar u,\nabla\bar P)$
is Lipschitz  continuous from  $\cM(\dot B^{n/p-1}_{p,1})\times\dot B^{n/p-1}_{p,1}$
to $E_p(T).$
 \end{theo}

The regularity given by Theorem \ref{th:lagrangian} ensures  that the map defined in \eqref{lag} 
is defined globally (see the Appendix).  Coming back to the Eulerian formulation,
this will enable us to get the following result\footnote{We here consider only the
case of small data to simplify the presentation.}:
\begin{theo}\label{th:eulerian}
Under the above assumptions,  System \eqref{eq:euler} has a unique global solution $(\rho, u,\nabla P)$ with 
$\rho\in L_\infty(\R_+;\cM(\dot B^{n/p-1}_{p,1}))$ and $(u,\nabla P)\in E_p.$
 \end{theo}
 Let us make a few comments concerning the above assumptions.
\begin{itemize}
\item
The condition $1\leq p<2n$ is a consequence of the product laws in Besov spaces. 
Let us emphasize  that any
space $L_\infty\cap B^{n/q-1}_{q,\infty}$ with $q$ satisfying
\begin{equation}\label{eq:q}\frac 1q>\frac 1n-\frac1p\quad\hbox{and}\quad
\frac1q\geq \frac1p-\frac1n,
\end{equation}
 embeds in $\cM(\dot B^{n/p-1}_{p,1}),$
a consequence of  basic continuity results for the 
paraproduct operator (see \cite{BCD}). Hence
the above statement  improves those of   \cite{Abidi,AP}
as regards the uniqueness. In particular, one may take the initial velocity
in a Besov space with a \emph{negative} index of regularity, 
so that a highly oscillating ``large'' velocity may give rise to a unique global solution. 
\item In contrast with the results of \cite{Abidi,AP,D1} it is not clear that 
the above statements may be generalized as so if the viscosity depends on the density. 
Recall that in this case, the diffusion term in the momentum equation of \eqref{eq:euler} reads
$\div(\mu(\rho)(\nabla u+{}^T\nabla u))$ where $\mu$ is a given suitably smooth nonnegative
function. One may easily extend the above  statements to this case  under the stronger condition 
that $(\mu(\rho_0)-\mu(1))$ is small in $\cM(\dot B^{n/p}_{p,1}).$
\end{itemize}

The space   $\cM(\dot B^{n/p-1}_{p,1})$ is in fact much larger
than $L_\infty\cap B^{n/q-1}_{q,\infty}$ with $q$ satisfying \eqref{eq:q} (see e.g. \cite{MS}, Chap. 4). 
It contains characteristic functions of $C^1$  bounded domains, whenever
$p>n-1$ (see the proof in Lemma \ref{l:jump}). 
Hence,  our result applies  to  mixture of fluids, which is of course of great 
physical interest. 
 In addition, given that the constructed velocity field $u$  is divergence-free and admits a $C^1$ flow $X$  (see again the appendix), we deduce the following
result which emphasizes  the range of Theorem \ref{th:lagrangian} (we just state the case of small velocities to simplify the presentation):
\begin{cor}\label{c:jump}
Assume that $u_0\in \dot B^{n/p-1}_{p,1}$ with $\div u_0=0$ and $n-1<p<2n.$
Let $\Omega_0$ be a bounded $C^1$ domain of $\R^n.$ 
There exist two constants $c$ (depending only on $p$ and $n$)
and $c'$ (depending only on $p,$ $n$ and $\Omega_0$) such that  if 
\begin{equation}
\|u_0\|_{\dot B^{n/p-1}_{p,1}}\leq c\quad\hbox{and}\quad
\rho_0=1+\sigma \chi_{\Omega_0}\ \hbox{ with } |\sigma|\leq c',
\end{equation}
then System \eqref{eq:euler} has a unique global solution 
$(\rho,u,\nabla P)$ with $\rho\in L_\infty(\R_+;\cM(\dot B^{n/p-1}_{p,1}))$
and $(u,\nabla P)\in E_p.$ 
In addition, for all time, 
\begin{equation}
\rho(t)=1+\sigma \chi_{\Omega_t}
\quad\hbox{where }\ \Omega_t=X_u(t,\Omega_0).\end{equation}
 Besides,  the measure and the $C^1$ regularity of $\d\Omega_t$ are preserved
for all time. 
\end{cor}

 Let us give the main ideas of the proof of existence. Obviously, it suffices to find a fixed point
for the map $\Theta:(q,v)\mapsto(\rho, u)$ where $(\rho,u)$ stands  
for the solution to the \emph{linear} system
\begin{equation}\label{eq:lineareuler}
\left\{\begin{array}{l}
\d_t\rho+v\cdot\nabla\rho=0,\\[1ex]
q(\d_tu+v\cdot\nabla u)-\mu\Delta u+\nabla P=0,\\[1ex]
\div u=0,\\[1ex]
u|_{|t=0}=u_0.
\end{array}\right.
\end{equation}
Although it is possible to prove uniform estimates in $L_\infty(\R_+;\cM(\dot B^{n/p-1}_{p,1}))\times E_p$ for $(\rho,u,\nabla P),$
we do not know how to get stability estimates in the same space, owing to the hyperbolic nature
of the density equation. As a consequence,  the contraction 
mapping theorem does not apply. 

In the present paper, we shall rather define  the solution 
of the above system in the \emph{Lagrangian coordinates corresponding to $v.$}
For such coordinates, the density is time-independent. 
So, given some reference vector field $\bar v$ and pressure field $\nabla \bar Q$
with $(\bar v,\nabla\bar Q)\in E_p,$  one may  define 
$(\bar u,\nabla\bar P)$ to be the solution 
 \emph{in  the Lagrangian coordinates $y=X_{v}^{-1}(t,x)$ 
pertaining to $v$} of the linear system \eqref{eq:lineareuler}
(that $X_v$ is a $C^1$-diffeomorphism over $\R^n$ is proved in the appendix).


Let us give more details. We  assume that $|DX_v|\equiv1$ and  set
$$
\bar\rho(t,y):=\rho(t,X_v(t,y)),\quad
\bar P(t,y):= P(t,X_v(t,y))\ \hbox{ and }\ 
\bar u(t,y)=u(t,X_v(t,y)),
$$
where $(\rho,u,\nabla P)$ stands for a solution to \eqref{eq:lineareuler}. Then
we have $\bar\rho(t,\cdot)\equiv\rho_0$ and  (see the proof in appendix)  
\begin{equation}\label{eq:linearlagrangian}
\begin{array}{l}
\rho_0\d_t\bar u-\mu\div (A_v{}^T\!A_v\nabla_y\bar u)
+{}^T\!A_v\nabla_y\bar P=0,\\[1ex]
\div_y(A_v\bar u)=0\end{array}
\end{equation}
with
\begin{equation}\label{eq:flowbar}
A_v=(D_yX_v)^{-1}\quad\hbox{and}\quad
X_v(t,y)=y+\int_0^t\bar v(\tau,y)\,d\tau.
\end{equation}

Solving this linear system \emph{globally} 
turns out to be possible under some smallness condition over $\bar v$
and $\rho_0-1.$ 
This will enable us to define a self-map $\Phi:(\bar v,\nabla\bar Q)\mapsto (\bar u,\nabla \bar P)$
on $E_p.$ 
Then it will be only a matter of checking  that  if the data $\rho_0$ and $u_0$ satisfy 
a suitable smallness condition then the map 
$\Phi$ fulfills the assumptions of the standard Banach fixed point theorem. 
The key to that will be  estimates for the Stokes system 
in $E_p$ (see Proposition \ref{p:stokeswhole})
and a ``magic'' algebraic relation involving the second equation of \eqref{eq:linearlagrangian}
(which has been used before in e.g. \cite{Mu,MZ} in a different context).


\section{The linear theory} \label{s:linear}

Our proof of existence for the linear System \eqref{eq:linearlagrangian}
will be based on the following a priori estimates for the Stokes system, the proof of which may be found in \cite{DM2}\footnote{Because homogeneous $0$-th order
multipliers act on any homogeneous Besov spaces $\dot B^s_{p,1},$ 
one may take any indices $s$ and exponents $p\in[1,\infty],$ in Proposition \ref{p:stokeswhole}.}:
\begin{prop}\label{p:stokeswhole}
 Let $p\in [1,\infty]$ and $ s\in\R.$ Let  $u_0\in \dot B^s_{p,1}(\R^n)$
 and  $f \in L_1(0,T;\dot B^s_{p,1}(\R^n)).$ 
 Let $g:[0,T]\times\R^n\rightarrow \R$ be such that 
$$\nabla g \in L_1(0,T;\dot B^{s}_{p,1}(\R^n)),\quad
\d_tg=\div R\  \mbox{~~with~~}\ R \in L_1(0,T;\dot B^s_{p,1}(\R^n))$$ 
and that  the compatibility condition   $g|_{t=0}=\div u_0$ on $\R^n$ is satisfied.
\medbreak
Then System 
\begin{equation}\label{w1}
 \begin{array}{lcr}
 \d_tu-\mu \Delta u +\nabla P=f \qquad & \mbox{in} & (0,T)\times\R^n\\[5pt]
\div u = g & \mbox{in} & (0,T)\times\R^n \\[5pt]
u|_{t=t_0} = u_0 &  \mbox{on} & \R^n
\end{array}
\end{equation}
 has a unique solution $(u,\nabla P)$ with 
$$
u\in\cC([0,T);\dot B^s_{p,1}(\R^n))\quad\hbox{and}\quad 
\d_tu,\nabla^2u,\nabla P\in L_1(0,T;\dot B^s_{p,1}(\R^n))
$$
and the following estimate is valid:
\begin{eqnarray}\label{w2}
\qquad\|u\|_{L_\infty(0,T;\dot B^s_{p,1}(\R^n))}+
 \|\d_tu,\mu\nabla^2 u ,\nabla P\|_{ L_1(0,T;\dot B^s_{p,1}(\R^n))}\hspace{5cm}\nonumber\\
 \hspace{2cm}
\leq C(\|f,\mu\nabla g,R\|_{ L_1(0,T;\dot B^s_{p,1}(\R^n))}+\|u_0\|_{\dot B^s_{p,1}(\R^n)})
\end{eqnarray}
where $C$ is an absolute constant with no dependence on $\mu$ and $T.$
\end{prop}

\medskip 

Granted with the above statement, we rewrite System \eqref{eq:linearlagrangian} as
\begin{equation}\label{eq:w3}
\!\!\!\begin{array}{l}
\d_t\bar u-\mu\Delta\bar u+\nabla\bar P=(1\!-\!\rho_0)\d_t\bar u
+\mu\div((A_v{}^T\!A_v-\Id)\nabla \bar u)+(\Id-{}^T\!A_v)\nabla \bar P\\[1ex]
\div\bar u=\div((\Id-A_v)\bar u)\\[1ex]
\bar u|_{t=0}=u_0.\end{array}
\end{equation}
We assume that the  vector-field $\bar v$ from which 
$A_v$ and $DX_v$ are defined satisfies
\begin{equation}\label{w4} 
\bar v\in\cC(\R_+;\dot B^{n/p-1}_{p,1}),\quad
 \d_t\bar v,\nabla^2\bar v\in L_1(\R_+;\dot B^{n/p-1}_{p,1}),
 \quad |DX_v|\equiv1
\end{equation}
and that, for a small enough constant $c,$ 
\begin{equation}\label{eq:smallv}
\int_0^\infty\|D\bar v\|_{\dot B^{n/p}_{p,1}}\,dt\leq c.
\end{equation}

Even though this system is linear, it cannot be solved directly 
by means of Proposition \ref{p:stokeswhole} for the right-hand side depends on the solution itself. So in order to prove the 
existence of $\bar u,$ we shall look for a fixed point of the map 
$$
\Psi:(\bar w,\nabla\bar Q)\longmapsto (\bar u,\nabla\bar P)
$$
where  $(\bar w,\nabla\bar Q)\in E_p$  and
 $(\bar u,\nabla\bar P)$ stands for the solution of
\begin{equation}\label{eq:systemw}\left\{\begin{array}{l}
\d_t\bar u-\mu\Delta\bar u+\nabla\bar P=f(\bar w,\nabla\bar Q),\\[1ex]
\div\bar u=g(\bar w),\\\bar u|_{t=0}=u_0.\end{array}\right.
\end{equation}
Above,  $g(\bar w):=\div((\Id-A_v)\bar w)$ and 
\begin{equation}\label{eq:f}
f(\bar w,\nabla\bar Q):=(1-\rho_0)\d_t\bar w
+\mu\div((A_v{}^T\!A-\Id)\nabla \bar w)+(\Id-{}^T\!A_v)\nabla \bar Q.
\end{equation}

We claim that,  if $\bar v$ satisfies \eqref{w4} and  the smallness condition \eqref{eq:smallv}, 
then for any $(\bar w,\nabla\bar Q)$ in the space $E_p$ with $1\leq p<2n,$ 
the above system has a unique solution $(\bar u,\nabla\bar P)$ in $E_p$
and that, in addition, the map $\Psi$ fulfills the required conditions for applying 
the contracting mapping theorem.

The existence of  $(\bar u,\nabla\bar P)$ will stem from Proposition \ref{p:stokeswhole} once
it has been checked that $f(\bar w,\nabla\bar Q)$ and $g(\bar w)$ fulfill the required conditions. 
As regards $g(\bar w),$ this stems from the following ``magic formula''

\begin{equation}\label{eq:magic}
g(\bar w)=\div((\Id-A_v)\bar w)=D\bar w:(\Id-A_v),
\end{equation}
a consequence of the fact that $|DX_v|\equiv1$ 
(see Corollary \ref{c:magic} in the Appendix). 

\subsubsection*{Bounds for $g(\bar w)$}

Let us first check that $g(\bar w)\in L_1(\R_+;\dot B^{n/p}_{p,1}).$ 
As $D\bar w\in L_1(\R_+;\dot B^{n/p}_{p,1})$ and as, according to \eqref{eq:U2},
$\Id-A_v$ is in $L_\infty(\R_+;\dot B^{n/p}_{p,1}),$ this is 
a consequence of the fact that $\dot B^{n/p}_{p,1}$ is a Banach algebra
and that 
$$
g(\bar w)=D\bar w:(\Id-A_v).
$$
In addition, we get 
\begin{equation} \label{eq:U7}
\|g(\bar w)\|_{L_1(\R_+;\dot B^{n/p}_{p,1})}\lesssim \|D\bar v\|_{L_1(\R_+;\dot B^{n/p}_{p,1})}
 \|D\bar w\|_{L_1(\R_+;\dot B^{n/p}_{p,1})}.
\end{equation}

Next, we see that
$\d_t(g(\bar w))=\div R^1(\bar w)+\div R^2(\bar w)$ with 
$$
R^1(\bar w):=(\Id-A_v)\d_t\bar w\quad\hbox{and}\quad R^2(\bar w):=-\d_tA_v\,\bar w.
$$
So, according to \eqref{eq:U1},\eqref{eq:U3} and because  the product operator maps 
$\dot B^{n/p}_{p,1}\times \dot B^{n/p-1}_{p,1}$ in $\dot B^{n/p-1}_{p,1}$
whenever  $p<2n,$ 
 we see that 
$R^1(\bar w)$ and $R^2(\bar w)$ belong to $L_1(\R_+;\dot B^{n/p-1}_{p,1})$ and that 
\begin{eqnarray}\label{eq:U5}
&&\|R^1(\bar w)\|_{L_1(\R_+;\dot B^{n/p-1}_{p,1})}\lesssim
 \|D\bar v\|_{L_1(\R_+;\dot B^{n/p}_{p,1})} \|\d_t\bar w\|_{L_1(\R_+;\dot B^{n/p-1}_{p,1})},\\
 \label{eq:U6}
 &&\|R^2(\bar w)\|_{L_1(\R_+;\dot B^{n/p-1}_{p,1})}\lesssim
 \|D\bar v\|_{L_1(\R_+;\dot B^{n/p}_{p,1})} \|\bar w\|_{L_\infty(\R_+;\dot B^{n/p-1}_{p,1})}.
\end{eqnarray}

\subsubsection*{Bounds for $f(\bar w,\nabla \bar Q)$}
 
That the first term of $f(\bar w)$  belongs to $L_1(\R_+;\dot B^{n/p-1}_{p,1})$
is a consequence of the definition  of the multiplier space
$\cM(\dot B^{n/p-1}_{p,1})$;
in addition we have
\begin{equation}\label{eq:U8}
\|(1-\rho_0)\d_t\bar w\|_{\dot B^{n/p-1}_{p,1}}\leq 
\|1-\rho_0\|_{\cM(\dot B^{n/p-1}_{p,1})}\|\d_t\bar w\|_{\dot B^{n/p-1}_{p,1}}.
\end{equation}

Next, according to \eqref{eq:U4}, the second term of $f(\bar w)$ belongs to 
$L_1(\R_+;\dot B^{n/p-1}_{p,1})$ and
\begin{equation}\label{eq:U9}
\|\div (A_v{}^T\!A_v-\Id)\nabla\bar w)\|_{L_1(\R_+;\dot B^{n/p\!-\!1}_{p,1})}\lesssim
 \|D\bar v\|_{L_1(\R_+;\dot B^{n/p}_{p,1})} \|D\bar w\|_{L_1(\R_+;\dot B^{n/p}_{p,1})}.
 \end{equation}
 
Finally, Inequality \eqref{eq:U2} 
and the fact that the product operator maps 
$\dot B^{n/p}_{p,1}\times \dot B^{n/p-1}_{p,1}$ in $\dot B^{n/p-1}_{p,1}$
if  $p<2n$ 
ensure that 
\begin{equation}\label{eq:U10}
\|(\Id-{}^T\!A_v)\nabla \bar Q\|_{L_1(\R_+;\dot B^{n/p-1}_{p,1})}\lesssim
 \|\nabla\bar v\|_{L_1(\R_+;\dot B^{n/p}_{p,1})} \|\nabla\bar Q\|_{L_1(\R_+;\dot B^{n/p-1}_{p,1})}.
\end{equation}

So putting together \eqref{eq:U7} to \eqref{eq:U10}, 
one may conclude from Proposition \ref{p:stokeswhole}
that  for any $(\bar w,\nabla\bar Q)$ in $E_p,$
System \eqref{eq:systemw} has a unique solution $(\bar u,\nabla\bar P)$ in $E_p.$ In addition
$$
\|(\bar u,\nabla\bar P)\|_{E_p}\leq C\Bigl(\|u_0\|_{\dot B^{n/p-1}_{p,1}}
+\bigl(\|1-\rho_0\|_{\cM(\dot B^{n/p-1}_{p,1})}
+\|D\bar v\|_{L_1(\R_+;\dot B^{n/p}_{p,1})}\bigr)
\|(\bar w,\nabla\bar Q)\|_{E_p}\Bigr).
$$

As a consequence, there exists a positive constant $c$ (depending only on $n$ and on $p$)
such that if \eqref{eq:smallv} is satisfied and 
\begin{equation}\label{eq:smallrho}
\|1-\rho_0\|_{\cM(\dot B^{n/p-1}_{p,1})}\leq c,
\end{equation}
then 
\begin{equation}\label{eq:U11}
\|\Psi(\bar w,\nabla\bar Q)\|_{E_p}\leq C\|u_0\|_{\dot B^{n/p-1}_{p,1}}
+\frac12\|(\bar w,\nabla\bar Q)\|_{E_p}.
\end{equation}
Banach theorem thus entails that the linear map $\Psi$ admits a unique fixed point in $E_p,$
that we shall still denote by $(\bar u,\nabla\bar P).$
Let us emphasize that  Inequality \eqref{eq:U11} ensures   that
\begin{equation}\label{eq:boundu}
\|(\bar u,\nabla\bar P)\|_{E_p}\leq 2C\|u_0\|_{\dot B^{n/p-1}_{p,1}},
\end{equation}
and that, by construction, $\div A_v\bar u=0$ hence, according to Corollary \ref{c:magic}, 
$|DX_u|\equiv1.$
In other words, given $\bar v$  fulfilling \eqref{w4} and \eqref{eq:smallv}, System \eqref{eq:w3} admits 
a unique solution $(\bar u,\nabla\bar P)$ fulfilling
 the same conditions, together with \eqref{eq:boundu}.


\section{The inhomogeneous Navier-Stokes equations}\label{s:nonlinear}

Let us denote by $\tilde E_{p}^R$ the closed subset of $E_p$ containing all the
couples  
$(\bar v,\nabla\bar Q)$ such that 
$$|DX_v|\equiv1\quad\hbox{and}\quad
\|(\bar v,\nabla\bar Q)\|_{E_p}\leq R.
$$

According to the previous section, if one takes $(\bar v,\nabla\bar Q)$ in $E_p$
 with $\bar v$ satisfying  $|DX_v|\equiv1$ and \eqref{eq:smallv}, then 
  \eqref{eq:linearlagrangian} admits a solution $(\bar u,\nabla\bar P)$
  in the same space, such that   
  $|DX_u|\equiv1.$  Let 
 $\Phi(\bar v,\nabla\bar Q)$ denote this solution\footnote{Of course, it is independent 
of $\nabla\bar Q.$ However, prescribing the pressure is needed
so as to define a map from a subset of $E_p$ to itself.}. 
We claim  that if $u_0$ is small 
enough with respect to $\mu$ in $\dot B^{n/p-1}_{p,1},$ if the density $\rho_0$
satisfies \eqref{eq:smallrho} and if $R$ is small enough, then $\Phi$
admits a unique fixed point in  $\tilde E_p^R,$
as  a consequence of the  contracting  
mapping  theorem.

\subsection{Stability of a small ball of $E_p$ by $\Phi$}

Assume that \eqref{eq:smallrho} is satisfied and let us take
$R=c\mu.$ Then  \eqref{eq:smallv} holds true whenever 
$(\bar v,\nabla\bar Q)$ is in $\tilde E_p^R.$ 
Therefore  $(\bar u,\nabla\bar P):=\Phi(\bar v,\nabla\bar Q)$ satisfies \eqref{eq:boundu}.
So it is clear that if 
\begin{equation}\label{eq:smallu0}
2C\|u_0\|_{\dot B^{n/p-1}_{p,1}}\leq c\mu,
\end{equation}
with $C$ as in \eqref{eq:boundu},
then $(\bar u,\nabla\bar P)$ is in $\tilde E_p^R,$ too.

\subsection{Contraction properties} 

In this part, we  show that under Conditions \eqref{eq:smallrho} and \eqref{eq:smallu0}
(with a greater constant $C$ and smaller constant $c$ if needed), 
  the map $\Phi:\tilde E_p^R\rightarrow\tilde E_p^R$ is $1/2$-Lipschitz.  

So we are given  $(\bar v_1,\nabla\bar Q_1)$ and $(\bar v_2,\nabla\bar Q_2)$ 
in $\tilde E_p^R,$ and denote
$$
(\bar u_1,\nabla\bar P_1):=\Phi(\bar v_1,\nabla\bar Q_1)\quad\hbox{and}\quad
(\bar u_2,\nabla\bar P_2):=\Phi(\bar v_2,\nabla\bar Q_2).$$
Let  $X_1$ and $X_2$ be the flows associated to $\bar v_1$ and $\bar v_2.$
Set $A_i=(DX_i)^{-1}$ for $i=1,2.$ The equations satisfied by 
$\du:=\bar u_2-\bar u_1$ and  $\nabla\dP:=\nabla \bar P_2-\nabla \bar P_1$ read 
$$\left\{\begin{array}{l}
\d_t\du-\mu\Delta\du+\nabla\dP=\df:=\df_1+\df_2+\df_3+\mu\div\df_4+\mu\div\df_5\\[1ex]
\div\du=\dg:=\div\bigl((\Id-A_2)\du+(A_1-A_2)\bar u_1\bigr)
\end{array}\right.
$$
with 
$$\displaylines{
\df_1:=(1-\rho_0)\d_t\du,\quad\df_2:=(\Id-{}^T\!A_2)\nabla\dP,\quad
\df_3:={}^T\!(A_1-A_2)\,\nabla \bar P_1,\cr
\df_4:=\bigl(A_2{}^T\!A_2-A_1{}^T\!A_1\bigr)\,\nabla \bar u_1,\quad
\df_5:=\bigl(A_2{}^T\!A_2-\Id\bigr)\nabla\du.}
$$

Once again, bounding $(\du,\nabla\dP)$ will stem from Proposition \ref{p:stokeswhole}, which 
ensures that, if $\d_t\dg=\div\dR$ then 
\begin{equation}\label{eq:dU}
\|(\du,\nabla\dP)\|_{E_p}\lesssim \|\df\|_{L_1(\R_+;\dot B^{n/p-1}_{p,1})}
+\mu\|\dg\|_{L_1(\R_+;\dot B^{n/p}_{p,1})}+
\|\dR\|_{L_1(\R_+;\dot B^{n/p-1}_{p,1})}.
\end{equation}
So we  have to estimate $\df_1,\df_2,\df_3$  and  $\df_4,\df_5$ 
in  $L_1(\R_+;\dot B^{n/p-1}_{p,1})$ and $L_1(\R_+;\dot B^{n/p}_{p,1}),$ respectively.
First, from the definition of the multiplier space $\cM(\dot B^{n/p-1}_{p,1}),$ we readily have
\begin{equation}\label{eq:df1}
\|\df_1\|_{L_1(\R_+;\dot B^{n/p-1}_{p,1})}\leq \|\rho_0-1\|_{\cM(\dot B^{n/p-1}_{p,1})}
\|\d_t\du\|_{L_1(\R_+;\dot B^{n/p-1}_{p,1})}.
\end{equation}
Next,  using Inequalities  \eqref{eq:U2},\eqref{eq:U4},  and product laws in Besov space yields
\begin{eqnarray}\label{eq:df2}
&&\|\df_2\|_{L_1(\R_+;\dot B^{n/p-1}_{p,1})}\lesssim
 \|D\bar v_2\|_{L_1(\R_+;\dot B^{n/p}_{p,1})} \|D\dP\|_{L_1(\R_+;\dot B^{n/p-1}_{p,1})},\\
\label{eq:df5}
&&\|\df_5\|_{L_1(\R_+;\dot B^{n/p}_{p,1})}\lesssim
 \|D\bar v_2\|_{L_1(\R_+;\dot B^{n/p}_{p,1})} \|D\du\|_{L_1(\R_+;\dot B^{n/p}_{p,1})}.
 \end{eqnarray}
Inequality \eqref{eq:dA}  ensures  that 
\begin{equation}\label{eq:df3}
\|\df_3\|_{L_1(\R_+;\dot B^{n/p}_{p,1})}\lesssim
 \|D\dv\|_{L_1(\R_+;\dot B^{n/p}_{p,1})}
 \|D\bar P_1\|_{L_1(\R_+;\dot B^{n/p-1}_{p,1})}
\end{equation}
whereas Inequalities \eqref{eq:dA}, \eqref{eq:dAdj} yield
\begin{equation}\label{eq:df4}
\|\df_4\|_{L_1(\R_+;\dot B^{n/p}_{p,1})}\lesssim
 \|D\dv\|_{L_1(\R_+;\dot B^{n/p}_{p,1})}
 \|D\bar u_1\|_{L_1(\R_+;\dot B^{n/p}_{p,1})}.
\end{equation}

In order to  bound $\dg$ in $L_1(\R_+;\dot B^{n/p}_{p,1}),$
we shall  use the fact that, by construction, 
$$\div \bar u_i=\div\bigl((\Id-A_i)\bar u_i\bigr)=D\bar u_i:(\Id-A_i).$$
Hence
$$
\dg=D\du:(\Id-A_2)-D\bar u_1:(A_2-A_1).
$$
  Now, easy  computations based on \eqref{eq:U2} and \eqref{eq:dA} yield
\begin{eqnarray}\label{eq:dg1}
&&\|D\du:(\Id-A_2)\|_{L_1(\R_+;\dot B^{n/p}_{p,1})}\lesssim
 \|D\bar v_2\|_{L_1(\R_+;\dot B^{n/p}_{p,1})} \|D\du\|_{L_1(\R_+;\dot B^{n/p}_{p,1})},\\\label{eq:dg2}
&&\|D\bar u_1:(A_2-A_1)\|_{L_1(\R_+;\dot B^{n/p}_{p,1})}\lesssim
 \|D\bar u_1\|_{L_1(\R_+;\dot B^{n/p}_{p,1})} \|D\dv\|_{L_1(\R_+;\dot B^{n/p}_{p,1})}.
\end{eqnarray}

Finally, to bound  $\d_t\dg,$ we decompose it into 
$\div(\dR_1+\dR_2+\dR_3+\dR_4)$ with 
$$
\begin{array}{lll}
&\dR_1=-\d_tA_2\:\du,\quad &\dR_2=(\Id-A_2)\d_t\du,\\[1ex]
&\dR_3=\d_t(A_1-A_2)\: \bar u_1,\quad
&\dR_4= (A_1-A_2)\d_t\bar u_1.\end{array}
$$
Using \eqref{eq:U1}, \eqref{eq:U3} and product laws in Besov spaces, we see that
\begin{eqnarray}\label{eq:dR1}
 &&\|\dR_1\|_{L_1(\R_+;\dot B^{n/p-1}_{p,1})}\lesssim
 \|D\bar v_2\|_{L_1(\R_+;\dot B^{n/p}_{p,1})}
 \|\du\|_{L_\infty(\R_+;\dot B^{n/p-1}_{p,1})},\\\label{eq:dR2}
  &&\|\dR_2\|_{L_1(\R_+;\dot B^{n/p-1}_{p,1})}\lesssim
  \|D\bar v_2\|_{L_1(\R_+;\dot B^{n/p}_{p,1})}
 \|\d_t\du\|_{L_1(\R_+;\dot B^{n/p-1}_{p,1})}.\end{eqnarray}
 
 In order to bound $\dR_3,$ it suffices to take advantage of \eqref{eq:dtdAdj}. 
 We get 
 \begin{equation}\label{eq:dR3}
   \|\dR_3\|_{L_1(\R_+;\dot B^{n/p-1}_{p,1})}\lesssim
   \|D\dv\|_{L_1(\R_+;\dot B^{n/p}_{p,1})}
\|\bar u_1\|_{L_\infty(\R_+;\dot B^{n/p-1}_{p,1})}.
\end{equation}
  Finally, using again \eqref{eq:dAdj}, we see that
    \begin{equation}\label{eq:dR4}
  \|\dR_4\|_{L_1(\R_+;\dot B^{n/p-1}_{p,1})}\lesssim  
  \|D\dv\|_{L_1(\R_+;\dot B^{n/p}_{p,1})}
  \|\d_t\bar u_1\|_{L_1(\R_+;\dot B^{n/p-1}_{p,1})}.
\end{equation}
 One can now plug Inequalities \eqref{eq:df1} to \eqref{eq:dR4} 
 in \eqref{eq:dU}. We end up with
 $$\displaylines{
 \|(\du,\nabla\dP)\|_{E_p}\leq C\bigl(\|1-\rho_0\|_{\cM(\dot B^{n/p-1}_{p,1})}
 +\|D\bar v_2\|_{L_1(\R_+;\dot B^{n/p}_{p,1})}\bigr) \|(\du,\nabla\dP)\|_{E_p}
 \hfill\cr\hfill+
 C\mu^{-1}\|(\bar u_1,\nabla\bar P_1)\|_{E_p}\|(\dv,\nabla\dQ)\|_{E_p}.}
 $$  
 So we see that if \eqref{eq:smallv} and \eqref{eq:smallrho} are  satisfied for 
 $\bar v_1,\bar v_2$ and $\rho_0$ with a small enough constant $c$, then we have
 $$
  \|(\du,\nabla\dP)\|_{E_p}\leq 2CR\mu^{-1}\|(\dv,\nabla\dQ)\|_{E_p}).
  $$
  Hence, the map $\Phi:\tilde E_p^R\mapsto\tilde E_p^R$ is $1/2$-Lipschitz whenever $R$
  and the data have been chosen so that \eqref{eq:smallrho}, \eqref{eq:smallu0}  are satisfied and 
  $4CR\leq\mu.$  
   This completes the proof of existence of a unique solution to System \eqref{eq:lagrangian} in 
   $\tilde E_p^R.$


\subsection{Stability estimates  in $E_p$}

In this part, we want to prove stability estimates in $E_p$ for the solutions 
to \eqref{eq:lagrangian}. This will ensure both uniqueness and  that the flow map is Lipschitz.

So we consider  two initial divergence-free velocity fields $u_{0,1}$ and $u_{0,2}$
in $\dot B^{n/p-1}_{p,1},$ and densities $(\rho_{0,1},\rho_{0,2})$ 
in $\cM(\dot B^{n/p-1}_{p,1})$ satisfying \eqref{eq:smallrho} and \eqref{eq:smallu0}. 
We want to  compare two solutions 
$(\bar u_1,\nabla \bar P_1)$ and $(\bar u_2,\nabla \bar P_2)$
in $E_p$   of   System \eqref{eq:lagrangian},  corresponding
to data $(\rho_{0,1},u_{0,1})$ and $(\rho_{0,2},u_{0,2}).$
 
The proof is similar to that of the contractivity of $\Phi$: we have to bound 
$$
\dU(t):= \|\du\|_{L_\infty(0,t;\dot B^{n/p-1}_{p,1})}+\|\d_t\du,\mu\nabla^2\du,\nabla\dPi\|_{L_1(0,t;\dot B^{n/p-1}_{p,1})}
$$
(with $\du:=\bar u_2-\bar u_1$ and  $\nabla\dP:=\nabla \bar P_2-\nabla \bar P_1$) 
in terms of $\|\du_{0}\|_{\dot B^{n/p-1}_{p,1}}$ and $\|\dr_{0}\|_{\cM(\dot B^{n/p-1}_{p,1})}.$
\smallbreak
Now,  the system for  $(\du,\nabla\dP)$ reads
$$\left\{\begin{array}{l}
\d_t\du-\mu\Delta\du+\nabla\dP=\df_0+\df_1+\df_2+\df_3+\mu\div(\df_4+\df_5),\\[1ex]
\div\du=\div\bigl((\Id-A_2)\du+(A_1\!-\!A_2)\bar u_1\bigr)
=D\du:(\Id-A_2)+D\bar u_1:(A_1\!-\!A_2)
\end{array}\right.
$$
with $\df_0:=\dr_0\,\d_t\bar u_1$ and where $\df_i$ for $i\in\{1,\cdots,5\}$
has been defined in the previous subsection. 
Of course, now the matrices $A_1$ and $A_2$ correspond to the vector-fields
$\bar u_1$ and $\bar u_2.$

\smallbreak
Using Proposition \ref{p:stokeswhole}, we gather that for all $t\in[0,T)$  
\begin{eqnarray}
\dU(t)\lesssim \|\du_{0}\|_{\dot B^{n/p\!-\!1}_{p,1}}
+ \sum_{i=0}^3 \|\df_i\|_{L_1(0,t;\dot B^{n/p\!-\!1}_{p,1})}
+\mu\sum_{i=4}^5 \|\df_i\|_{L_1(0,t;\dot B^{n/p}_{p,1})}
+\! \sum_{i=1}^4 \|\dR_i\|_{L_1(0,t;\dot B^{n/p\!-\!1}_{p,1})}\quad\nonumber\\\label{eq:dU1}\hspace{4.4cm}
+\mu\|D\du:(\Id\!-\!A_2)\|_{L_1(0,t;\dot B^{n/p}_{p,1})}
+\mu\|D\bar u_1:(A_2\!-\!A_1)\|_{L_1(0,t;\dot B^{n/p}_{p,1})}.
\end{eqnarray}
{}From the definition of the multiplier space $\cM(\dot B^{n/p-1}_{p,1}),$ we readily have
$$
\|\df_0\|_{\dot B^{n/p-1}_{p,1}}\leq\|\dr_0\|_{\cM(\dot B^{n/p-1}_{p,1})}\|\d_t\bar u_1\|_{\dot B^{n/p-1}_{p,1}}.$$
The other terms may be bounded as in the previous subsection. 
So we eventually conclude that for all $t\in[0,T)$
 $$
 \dU(t)\leq \frac12\dU(t)+C\bigl(\|\dr_0\|_{\cM(\dot B^{n/p-1}_{p,1})}+\|\du_0\|_{\dot B^{n/p-1}_{p,1}}\bigr)
 $$  
  whenever, for $i=1,2,$ 
  $$
  \sup_{t\in[0,T]} \|\bar u_i(t)\|_{\dot B^{n/p-1}_{p,1}}+
  \int_0^T\Bigl(\mu^{-1}\|\d_t\bar u_i,\nabla\bar P_i\|_{\dot B^{n/p-1}_{p,1}}
  +\|D\bar u_i\|_{\dot B^{n/p}_{p,1}}\Bigr)\,dt
  $$
  is small enough. 
  \medbreak
  This latter condition is a consequence of \eqref{eq:smallu0}. 
  This completes the proof of stability estimates.


  \subsection{Proof of the local-in-time existence result}
  
  We here explain how the arguments of the previous subsections have to be modified so as to handle 
  large initial velocities. 
  
   Let us first notice that the computations that have been performed in Section \ref{s:linear} 
  also hold \emph{locally} on $[0,T)$ whenever $\bar v$ 
  satisfies
  \begin{equation}\label{w4T}
  \bar v\in \cC_b([0,T);\dot B^{n/p-1}_{p,1}),\quad
  \d_t\bar v,\nabla^2\bar v\in L_1(0,T;\dot B^{n/p-1}_{p,1}),\quad
  |DX_v|\equiv1\ \hbox{ on }\ [0,T)\times\R^n
  \end{equation}
  and
  \begin{equation}\label{eq:smallvT}
  \int_0^T\|D\bar v\|_{\dot B^{n/p}_{p,1}}\,dt\leq c.
  \end{equation} 
  This ensures that, under Condition \eqref{eq:smallrho},
   System \eqref{eq:linearlagrangian} may be solved locally in $E_p(T).$ Of course, Inequality
   \eqref{eq:boundu} is still satisfied (for the norm in $E_p(T)$). However, it is not accurate
   enough so as to solve the nonlinear system, if $u_0$ is too large. 
   To overcome this, we shall apply the contracting mapping
  theorem  in some suitable neighborhood of the solution $(u_L,\nabla P_L)$ to the ``free'' Stokes system, that is
  \begin{equation}\label{eq:free}
  \begin{array}{lcr}
  \d_tu_L-\mu\Delta u_L+\nabla P_L=0&\hbox{ in }& [0,T)\times\R^n,\\[1ex]
  \div u_L=0&\hbox{ in }& [0,T)\times\R^n,\\[1ex]
  u_L|_{t=0}=u_0&\hbox{ on }&\R^n.
  \end{array}
  \end{equation}
  Setting $(\bar u,\nabla\bar P):=\Phi(\bar v,\nabla\bar Q)$, we want to show that
  if $T$ is small enough (a condition which will be expressed in terms of
  the free solution only) then $(\tilde u,\nabla\tilde P):=(\bar u-u_L,\nabla(\bar P-P_L))$ is small. 
  For that, we shall apply Proposition  \ref{p:stokeswhole} to the system satisfied
  by $(\tilde u,\nabla\tilde P),$ namely 
  \begin{equation}
\!\!\!\begin{array}{l}
\d_t\tilde u-\mu\Delta\tilde u+\nabla\tilde P=f(\bar u,\nabla \bar P)
\\[1ex]
\div\tilde u=g(\bar u),\\[1ex]
\tilde  u|_{t=0}=0\end{array}
\end{equation}
  where $f(\bar u,\nabla\bar P)$ and $g(\bar u)$ have been defined in Section \ref{s:linear}.
  \smallbreak
  On the one hand, we shall bound $f(\bar u,\nabla\bar P)$ in $L_1(0,T;\dot B^{n/p-1}_{p,1})$
  and $g(\bar u)$ in $L_1(0,T;\dot B^{n/p}_{p,1})$ exactly as in Section \ref{s:linear}, 
  on the other hand, decomposing $\d_t(g(\bar u))$
  into $R^1(\bar u)+R^2(\bar u),$ we see that the bound \eqref{eq:U6} 
  for $R^2(\bar u)$ is not accurate enough as it involves 
  $\|\bar u\|_{L_\infty(0,T;\dot B^{n/p-1}_{p,1})}$ which need not be small for $T$ going
  to $0,$ if $u_0$ is large. 
  So we shall rather write
  $$ 
  R^2(\bar u)=-\d_tA_v \:u_L-\d_tA_v\:\tilde u,
  $$
  and use product laws and Inequality \eqref{eq:U3b}, to get
  $$
  \|R^2(\bar u)\|_{L_1(0,T;\dot B^{n/p-1}_{p,1})}\lesssim
  \|D\bar v\|_{L_2(0,T;\dot B^{n/p-1}_{p,1})}\|u_L\|_{L_2(0,T;\dot B^{n/p}_{p,1})}
+ \|D\bar v\|_{L_1(0,T;\dot B^{n/p}_{p,1})}\|\tilde u\|_{L_\infty(0,T;\dot B^{n/p-1}_{p,1})}.
$$
  
  So, finally, using Proposition  \ref{p:stokeswhole} and decomposing
  everywhere $\bar u$ and $\nabla\bar P$ in $u_L+\tilde u$
  and $\nabla P_L+\nabla \tilde P,$ we get
  $$\displaylines{
  \|(\tilde u,D\tilde P)\|_{E_p(T)}\leq C
  \bigl(\|1-\rho_0\|_{\cM(\dot B^{n/p-1}_{p,1})} +\|D\bar v\|_{L_1(0,T;\dot B^{n/p}_{p,1})}\bigr)
   \hfill\cr\hfill\times \bigl(\|(\tilde u,D\tilde P)\|_{E_p(T)}
    +\|\d_tu_L,\mu D^2u_L,DP_L\|_{L_1(0,T;\dot B^{n/p-1}_{p,1})}\bigr)
  +C \|\bar v\|_{L_2(0,T;\dot B^{n/p}_{p,1})}\|u_L\|_{L_2(0,T;\dot B^{n/p}_{p,1})}.}
   $$
  Therefore, if \eqref{eq:smallrho} and \eqref{eq:smallvT} are satisfied with $c$
  small enough and denote $\tilde v:=\bar v-u_L,$ then we get 
  $$
  \displaylines{
   \|(\tilde u,D\tilde P)\|_{E_p(T)}\leq 
   Cc\|\d_tu_L,\mu D^2u_L,DP_L\|_{L_1(0,T;\dot B^{n/p-1}_{p,1})}\hfill\cr\hfill
   +C\|u_L\|_{L_2(0,T;\dot B^{n/p}_{p,1})}^2
    +C \|\tilde v\|_{L_2(0,T;\dot B^{n/p}_{p,1})}\|u_L\|_{L_2(0,T;\dot B^{n/p}_{p,1})}.}
  $$
  This inequality together with the interpolation inequality
  $$ 
  \|\tilde v\|_{L_2(0,T;\dot B^{n/p}_{p,1})}\|\leq \|\tilde v\|_{L_1(0,T;\dot B^{n/p+1}_{p,1})}^{1/2}
   \|\tilde v\|_{L_\infty(0,T;\dot B^{n/p-1}_{p,1})}^{1/2}
  $$
  ensures that $\Phi$ maps $(u_L,\nabla P_L)+\tilde E_p^R(T)$
  (where $\tilde E_p^R(T)$ is the ``local'' version of $\tilde E_p^R$) into itself whenever
  $T$ satisfies
  \begin{equation}\label{eq:smallfree}
  \begin{array}{l}
   Cc\|\d_tu_L,\mu D^2u_L,DP_L\|_{L_1(0,T;\dot B^{n/p-1}_{p,1})}
      +C\|u_L\|_{L_2(0,T;\dot B^{n/p}_{p,1})}^2\leq R/2,\\[1ex]
   C\mu^{-1/2}\|u_L\|_{L_2(0,T;\dot B^{n/p}_{p,1})}\leq1/2.\end{array}
   \end{equation}
  Of course, for  \eqref{eq:smallvT} to be satisfied, it suffices that 
  to take $R=c\mu/2$ and to assume that 
  $T$ is so small as
  \begin{equation}\label{eq:smallfreebis}
  \|Du_L\|_{L_1(0,T;\dot B^{n/p}_{p,1})}\leq c/2.
  \end{equation}
  So if   \eqref{eq:smallfree} and \eqref{eq:smallfreebis} 
  are satisfied (conditions which depend only on the data) then  one may conclude
  that $\Phi$ maps $(u_L,\nabla P_L)+\tilde E_p^R(T)$ into itself. 
  \smallbreak
  The proof of the contraction properties for $\Phi$ in this context follows the same lines : 
  we consider $(\bar u_i,\nabla\bar P_i)=\Phi(\bar v_i,\nabla\bar Q_i)$
  with $(\bar v_i,\nabla \bar Q_i)$ in $(u_L,\nabla P_L)+\tilde E_p^R(T)$ for $i=1,2,$ then  
  we bound all the terms $\df_i,$ $\dg_i$ and $\dR_i$ as in the case of
  small initial velocity, except for $\dR_3$ as it involves
 $\|\bar u_1\|_{L_\infty(0,T;\dot B^{n/p-1}_{p,1})}$ which need not be small for $T$ going to $0.$
  For this latter term, we notice that, according to \eqref{eq:dtdAdjb},
 $$
 \|\dR_3\|_{L_1(0,T;\dot B^{n/p-1}_{p,1})}\lesssim 
\|\dv\|_{L_2(0,T;\dot B^{n/p}_{p,1})} \|\bar u_1\|_{L_2(0,T;\dot B^{n/p}_{p,1})}.
$$
So we eventually get, 
 $$\displaylines{
 \|(\du,\nabla\dP)\|_{E_p(T)}\leq C\bigl(\|1-\rho_0\|_{\cM(\dot B^{n/p-1}_{p,1})}
 +\|D\bar v_2\|_{L_1(0,T;\dot B^{n/p}_{p,1})}\bigr) \|(\du,\nabla\dP)\|_{E_p(T)}
 \hfill\cr\hfill+
 C\bigl(\mu^{-1}\|(\tilde u_1,\nabla\tilde P_1)\|_{E_p(T)}
 +\mu^{-1/2}\|u_L\|_{L_2(0,T;\dot B^{n/p}_{p,1})}\bigr)
 \|(\dv,\nabla\dQ)\|_{E_p(T)}).}
 $$  
Note that our assumptions on $(\tilde u_1,\nabla\tilde P_1)$ and on the free solution ensures that
if $R$ and $T$ have been chosen small enough then the factor of the last term
is smaller than, say, $1/2.$ So the contraction mapping theorem applies. 
This completes the proof of the existence part of Theorem \ref{th:local}. 
Proving the  stability and uniqueness follows from similar arguments. 
The details are left to the reader.


  \subsection{Proof of Theorem \ref{th:eulerian}}

  Given data $(\rho_0,u_0)$ satisfying the assumptions of Theorem 
  \ref{th:eulerian}, one may construct a global solution $(\bar u,\nabla\bar P)$
  to System \eqref{eq:lagrangian} in $E_p.$ 
  If $X_u$ denotes the ``flow'' to $\bar u$ which is defined according 
  to \eqref{eq:flowbar} 
  then the results of the appendix ensure that, for all $t\in\R_+,$
  $X_u(t,\cdot)$ is a $C^1$ diffeomorphism of $\R^n.$
  In particular, one may set 
  $$
  \rho(t,\cdot):=\rho_0\circ X^{-1}_u(t,\cdot),\quad
  P(t,\cdot):=\bar P(t,\cdot)\circ  X^{-1}_u(t,\cdot),\quad
  u(t,\cdot):=\bar u(t,\cdot) \circ X^{-1}_u(t,\cdot),
  $$
  and the algebraic relations that are derived in the 
  appendix show  that $(\rho,u,\nabla P)$ satisfies 
 System \eqref{eq:euler}. 
 In addition, given that $X_u(t,\cdot)$ is measure
 preserving and that $DX_u(t)-\Id$ belongs to $\dot B^{n/p}_{p,1},$
 the map $a\mapsto a\circ X^{\pm1}_u(t)$ is continuous from 
  $\dot B^{s}_{p,1}$ to itself if  $s\in\{n/p-1,n/p\},$
 (see e.g. \cite{DM3}, Chap. 2). 
 This implies that:
 \begin{itemize}
 \item the Eulerian velocity $u$ is in $\cC_b(\R_+;\dot B^{n/p-1}_{p,1}),$
 \item for any $\phi\in \dot B^{n/p-1}_{p,1},$
 we have $\phi\circ X_u^{\pm1}(t)\in \dot B^{n/p-1}_{p,1}.$
    So, given that $\rho_0\in \cM(\dot B^{n/p-1}_{p,1}),$ we have
  $$
  \phi\,\rho(t)=\bigl((\phi\circ X_u(t))\rho_0\bigr)\circ X_u^{-1}(t) \in \dot B^{n/p-1}_{p,1}.
  $$
  Hence $\rho\in L_\infty(\R_+;\cM(\dot B^{n/p-1}_{p,1})).$
  \item the chain rule ensures that
  $$
  \nabla P=({}^T\! A_u\cdot\nabla\bar P)\circ X_u^{-1}.
  $$
  So combining product laws and the invariance of $\dot B^{n/p-1}_{p,1}$
  by right-composition, we get
  $\nabla P\in L_1(\R_+;\dot B^{n/p-1}_{p,1}).$
  \item the chain rule also ensures that
  $$
  \nabla u=({}^T\! A_u\cdot\nabla\bar u)\circ X_u^{-1}.
  $$
  So using the fact that $\dot B^{n/p}_{p,1}$ is a Banach algebra, 
  and using invariance of $\dot B^{n/p}_{p,1}$
  by right-composition, we get
  $\nabla u\in L_1(\R_+;\dot B^{n/p}_{p,1}).$
\end{itemize}  

In order to prove uniqueness, we consider two solutions
$(\rho_1,u_1,\nabla P_1)$ and $(\rho_2,u_2,\nabla P_2)$ corresponding
to the same data $(\rho_0,u_0),$ and 
perform the Lagrangian  change of variable (pertaining to the flow of 
$u_1$ and $u_2$ respectively. 
The obtained functions $(\bar u_1,\nabla\bar P_1)$ and 
$(\bar u_2,\nabla\bar P_2)$ both satisfy \eqref{eq:lagrangian}
with the same $\rho_0$ and $u_0.$ Hence they coincide, as 
a consequence of the uniqueness part of  Theorem \ref{th:lagrangian}.


\appendix
\section{Appendix}
\setcounter{equation}{0}

Let us first derive  algebraic relations involving changes of  coordinates. 

We are  given a  $C^1$-diffeomorphism $X$ over $\R^n.$
For $H:\R^n\rightarrow\R^m,$ we agree that  $\bar H(y)= H(x)$ with $x=X(y).$ 
With this convention,  the chain rule writes
\begin{equation}\label{eq:chainrule}
D_y\bar H(y)=D_xH(X(y))\cdot D_yX(y)
\quad\hbox{with }\  (D_yX)_{ij}=\d_{y_j}X^i.
\end{equation}
or, denoting $\nabla_y={}^T\!D_y,$ 
$$
\nabla_y\bar H(y)= (\nabla_yX(y))\cdot\nabla_xH(X(y)).
$$

Hence  we have 
\begin{equation}\label{eq:A}
D_xH(x)=D_y\bar H(y)\cdot A(y)\quad\hbox{with}\quad A(y)=(D_yX(y))^{-1}=D_xX^{-1}(x).
\end{equation}

\begin{lem}\label{l:div} Let  $H$ be a vector-field over $\R^n.$ 
If we denote $\bar H=H\circ X$ then the following relation holds true: 
\begin{equation}\label{eq:div}
\div_x  H(x)=\div_y(A|DX|\bar H)(y)=\div_y(\adj(DX)\bar H)(y)\quad\hbox{with }\ x=X(y),
\end{equation}
where $|DX|$ stands for the determinant of $DX,$ and  $\adj(DX)$ for the
adjugate of $DX,$ that is the transpose of the cofactor matrix of $DX.$
\end{lem}
\begin{p}
This stems  from the following series of computations (based on integrations by parts
and \eqref{eq:A})  which hold for any
scalar test function $q$:  
$$
\begin{array}{lll}
\Int q(x)\div_x H(x)\,dx&=&-\Int D_xq(x)\cdot H(x)\,dx,\\[1ex]
&=&-\Int D_xq(X(y))\cdot H(X(y))|D_yX(y)|\,dy\\[1ex]
&=&-\Int D_y\bar q(y)\cdot (A|D_yX|)(y)\cdot \bar H(y)\,dy,\\[1ex]
&=&\Int \bar q(y)\div_y(A|D_yX| H)(y)\,dy.\end{array}
$$
As $A=|D_yX|^{-1}\adj(D_yX),$ we get the result.
\end{p}
\begin{rem}
Combining \eqref{eq:A} and \eqref{eq:div}, we deduce that  if $a:\R^n\rightarrow\R$ then
$$
\overline{\Delta_x a}=\overline{\div_x\nabla_x a}=\div_y(A|DX|\overline{\nabla_xa})
=\div_y (\adj(DX){}^T\!A\nabla_y\bar a).
$$
\end{rem}

Recall that if  $v$ is a  time-dependent vector field with coefficients
in $L_1(0,T;C^{0,1})$ then it has, by virtue of the Cauchy-Lipschitz theorem, a unique $C^1$
flow $X_v$ satisfying
$$
X_v(t,y)=y+\int_0^t v(\tau,X_v(\tau,y))\,d\tau\quad\hbox{for all }\ t\in[0,T),
$$
 and that $X_v(t,\cdot)$ is a $C^1$-diffeomorphism over $\R^n.$ 
\medbreak
Lemma \ref{l:div}  enables us to deduce the following ``magic''  relation
which is the corner stone of the proof of our main results:
\begin{cor}\label{c:magic} 
Let $v$ and $w$ be two time-dependent vector fields with coefficients
in $L_1(0,T;C^{0,1}).$  Let $X_v$ and $X_w$ be the corresponding
flows.  Denote $A_v:=(DX_v)^{-1}$ and $A_w:=(DX_w)^{-1}.$ 
Let us introduce the Lagrangian coordinates $y_v$ and $y_w$ pertaining to $v$
and $w,$ respectively, defined by
$$
x=X_v(y_v)=X_w(y_w).
$$
Assume in addition that  $$|DX_v|\equiv1\ \hbox{ and }\   
\div(A_v\bar w_v)=0\quad\hbox{with}\quad
\bar w_v:=w\circ X_v.$$

Then $|DX_w|\equiv1$ and for any $C^1$   vector-field $H,$ one has
$$\div H(x)=(D \bar H_v:A_v)(y_v)=\div\bigl(A_w\bar H_w\bigr)(y_w)
\quad\hbox{with }\  \bar H_v:=H\circ X_v\ \hbox{ and }\ \bar H_w:=H\circ X_w.
$$
\end{cor}

\begin{p}
With the above notation, the chain rule ensures that 
$$
D_xH(x)=D_{y_v}\bar H_{v}(y_v)\cdot A_v(y_v).
$$
Hence taking the trace yields the left equality. 
\smallbreak
Next, according to Lemma \ref{l:div} and to our assumption over $v$ and $w,$ we have
$$
0=\div (A_v\bar w_v)=\div(\adj(DX_v)\bar w_v)=\div_x w(x).
$$
Hence Liouville theorem ensures that $|DX_w|\equiv1.$
So finally, applying  Lemma \ref{l:div} with $X_w$ completes the proof. 
\end{p}


\begin{lem}\label{l:adj} There exist
$n^2$ at least quadratic polynomials $P_{ij}:\cM_n(\R)\rightarrow\R$ of degree  $n-1$ such that
$$
\Id-\adj(\Id+C)=\bigl(C-({\rm Tr}\, C)\Id\bigr)+P_2(C),
$$
where $P_2(C)$ is the $n\times n$ matrix with entries $P_{ij}(C).$
\end{lem}
\begin{p}
It suffices to use the fact that,  by definition of the differential of $\adj,$ we have
$$
\Id-\adj(\Id+C)=\adj(\Id)-\adj(\Id+C)=-d\,\adj(\Id)(C)+P_2(C).
$$
Now, 
$$
\adj(\Id+C)=
(\Id+C)^{-1}\det(\Id+C)
$$
and the differential of the reciprocal operator
at $\Id$ is $C\mapsto -C$ while  $d\det(\Id)(C)=({\rm Tr}\, C)\Id.$
So $d\,\adj(\Id)(C)=({\rm Tr}\, C)\,\Id-C.$
\end{p}
\medbreak

We now want to establish some a priori estimates for the flow
which will be needed in our main results. 
The first difficulty that has to be faced is that when implementing
the iterative process for solving \eqref{eq:lagrangian}, we are given 
the velocity field $\bar v$ \emph{in Lagrangian coordinates}. 
Therefore, it first has to be checked whether the ``flow'' $X_v(t,\cdot)$ 
defined by 
\begin{equation}\label{eq:defX}
X_v(t,y):=y+\int_0^t\bar v(\tau,y)\,d\tau
\end{equation}
is a $C^1$ diffeomorphism over $\R^n.$
This property is required for constructing
the Eulerian vector-field $v$ by setting $v(t,\cdot):= v\circ X_v^{-1}(t,\cdot).$ 
\medbreak
So let us assume that we are given
 some vector field $\bar v$ over $[0,T)\times\R^n$ with 
$$
\bar v\in \cC_b([0,T);\dot B^{n/p-1}_{p,1}),\quad
\d_t\bar v\in L_1([0,T);\dot B^{n/p-1}_{p,1})\ \hbox{ and }\ 
D\bar v\in L_1([0,T);\dot B^{n/p}_{p,1}).
$$
Differentiating \eqref{eq:defX} with respect to the space variable yields
\begin{equation}\label{eq:DX}
DX_v(t,y):=\Id+\int_0^tD\bar v(\tau,y)\,d\tau.
\end{equation}

As $\dot B^{n/p}_{p,1}(\R^n)$ is embedded in the set $\cC_0(\R^n)$ 
of continuous functions going to $0$ at infinity, 
we deduce that  $X_v$  is a  $C^1$ function over
 $\R_+\times\R^n.$
 However, in general, $X_v(t,\cdot)$ need not be
  a $C^1$-diffeomorphism over $\R^n$ for all $t\in\R_+.$
 So we assume that the smallness condition \eqref{eq:smallv} is satisfied
 with $c$ small enough. Then, using embedding  we see that 
  it  guarantees that  
  $$\|DX_v(t,\cdot)-\Id\|_{L_\infty(\R^n)}\leq1/2\quad\hbox{for all }\ t\in\R_+.
$$
Hence, for any $t\in\R_+,$ the map  $X_v(t,\cdot)$ is a local 
diffeomorphism.
In order to show that it is a \emph{global} diffeomorphism, 
we introduce the solution $Y_v$ to   the ordinary differential equation
\begin{equation}\label{eq:ODE}
\bar v(t,Y_v(t,x))+DX_v(t,Y_v(t,x))\frac{d}{dt} Y_v(t,x)=0.
\end{equation}

Under \eqref{eq:smallv}, the matrix $DX_v$ is invertible at every point and
$(DX_v)^{-1}-\Id$ belongs to $L_\infty(\R_+;\dot B^{n/p}_{p,1}).$
Indeed, one may write
$$
(DX_v)^{-1}-\Id=\sum_{k\geq1}(\Id-DX_v)^k.
$$
Hence, using the assumptions over $\bar v$
and the product laws in Besov spaces (here we need  that $1\leq p<2n$),
  \eqref{eq:ODE} may be seen as an ordinary differential equation in the Banach space  $\dot B^{n/p-1}_{p,1},$
which may be solved on $[0,T)$ according to Cauchy-Lipschitz theorem. 
  In particular, differentiating  $X_v(t,\cdot)\circ Y_v(t,\cdot)$
with respect to time, we easily gather that  $X_v(t,\cdot)\circ Y_v(t,\cdot)=\Id$ for all $t\in[0,T).$ Therefore, one may eventually conclude that $X_v(t,\cdot)$ is a $C^1$-diffeomorphism
over $\R^n,$ with inverse $Y_v(t,\cdot).$
\smallbreak

Let us now derive some ``flow estimates'' that will be needed for constructing
the maps $\Phi$ and $\Psi.$ For completeness, the statements that we here state
are slightly more general than needed : that $|DX|\equiv1$ is not assumed.  
\begin{lem}\label{eq:flow} Let $p\in[1,+\infty).$ Under Assumption \eqref{eq:smallv} for $\bar v,$
we have
\begin{eqnarray}\label{eq:U1}
&&\|\Id-\adj(DX(t))\|_{\dot  B^{n/p}_{p,1}}\lesssim \|D\bar v\|_{L_t^1(\dot B^{n/p}_{p,1})},\\
\label{eq:U2}
&&\|\Id-A(t)\|_{\dot B^{n/p}_{p,1}}\lesssim \|D\bar v\|_{L_1(0,t;\dot B^{n/p}_{p,1})},\\\label{eq:U3}
&&\|\d_t(\adj (DX))(t)\|_{\dot B^{n/p}_{p,1}}\lesssim  \|D\bar v(t)\|_{\dot B^{n/p}_{p,1}},\\\label{eq:U3b}
&&\|\d_t(\adj (DX))(t)\|_{\dot B^{n/p-1}_{p,1}}\lesssim  \|D\bar v(t)\|_{\dot B^{n/p-1}_{p,1}}
\quad\hbox{if }\ p<2n,\\\label{eq:U4}
&&\|\adj(DX(t)){}^T\!A(t)-\Id\|_{\dot B^{n/p}_{p,1}}
\lesssim \|D\bar v\|_{L_1(0,t;\dot B^{n/p}_{p,1})}.
\end{eqnarray}
\end{lem}
\begin{p}
According to Lemma \ref{l:adj} and to \eqref{eq:DX}, one may write
$$
\Id-\adj(DX(t))=\int_0^t\bigl(D\bar v-\div\bar v\,\Id\bigr)\,d\tau
+P_2\biggl(\Bigl(\int_0^t D\bar v\,d\tau\Bigr)\biggr)
$$
where the coefficients of $P_2$ are at least quadratic polynomials of degree $n-1.$
Given that $\dot B^{n/p}_{p,1}$ is a Banach algebra and that \eqref{eq:smallv} holds,  we readily 
get the result. 

In order to prove the second estimate, we just use the fact that, under assumption \eqref{eq:smallv}, we have  
\begin{equation}\label{eq:C}
A(t)=(\Id+C(t))^{-1}=\sum_{k\in\N} (-1)^k(C(t))^k\quad\hbox{with}\quad
C(t)=\int_0^t D\bar v\,d\tau,
\end{equation}
and that $\dot B^{n/p}_{p,1}$ is a Banach algebra. 
\smallbreak
In order to prove the third inequality, we use the fact that, according to Lemma \ref{l:adj}, 
we have
$$
\d_t\bigl(\adj(DX)\bigr)=\frac\d{\d t}\biggl(\Id+\int_0^t\bigl(\div\bar v\,\Id-D\bar v\bigr)\,d\tau
+P_2\Bigl(\int_0^tD\bar v\,d\tau\Bigr)\biggr)\cdotp
$$
Hence
$$
\d_t\bigl(\adj(DX)\bigr)(t)=\bigl(\div\bar v(t)\Id-D\bar v(t)\bigr)+dP_2\Bigl(\int_0^tD\bar v\,d\tau\Bigr)\cdot D\bar v(t).
$$
As the coefficients of $dP_2$ are polynomials of $n^2$ variables   that vanish 
at $0,$ we get 
$$
\|\d_t(\adj (DX))(t)\|_{\dot B^{n/p}_{p,1}}\lesssim  
\|D\bar v(t)\|_{\dot B^{n/p}_{p,1}}\biggl(1+\biggl\|\int_0^tD\bar v\,d\tau\biggr\|_{\dot B^{n/p}_{p,1}}\biggr),
$$
hence \eqref{eq:U3}.
Proving \eqref{eq:U3b} is similar: it is only a matter of using the continuity of
the product from $\dot B^{n/p-1}_{p,1}\times\dot B^{n/p}_{p,1}$
to  $\dot B^{n/p-1}_{p,1},$ if $p<2n.$
\smallbreak

For proving the last inequality, we use the decomposition
$$
\adj(DX){}^T\!A-\Id=(\adj(DX)-\Id){}^T\!A+{}^T\!(A-\Id).
$$
So combining Inequalities \eqref{eq:U1} and \eqref{eq:U2}, and
the fact that $\dot B^{n/p}_{p,1}$ is a Banach algebra, we get
the result.
\end{p}

\begin{lem} 
Let $\bar v_1$ and $\bar v_2$ be two vector-fields satisfying \eqref{eq:smallv},
and $\dv:=\bar v_2-\bar v_1.$ 
Then we have for all  $p\in[1,+\infty),$
\begin{equation}\label{eq:dA}
\|A_2-A_1\|_{L_\infty(\R_+;\dot B^{n/p}_{p,1})} \lesssim 
 \|D\dv\|_{L_1(\R_+;\dot B^{n/p}_{p,1})},
\end{equation}
\begin{equation}\label{eq:dAdj}
\|\adj(DX_2)-\adj(DX_1)\|_{L_\infty(\R_+;\dot B^{n/p}_{p,1})} \lesssim 
 \|D\dv\|_{L_1(\R_+;\dot B^{n/p}_{p,1})},
\end{equation}
\begin{equation}\label{eq:dtdAdj}
\|\d_t(\adj(DX_2)-\adj(DX_1))\|_{L_1(\R_+;\dot B^{n/p}_{p,1})} \lesssim 
 \|D\dv\|_{L_1(\R_+;\dot B^{n/p}_{p,1})},
\end{equation}
\begin{equation}\label{eq:dtdAdjb}
\|\d_t(\adj(DX_2)-\adj(DX_1))\|_{L_2(\R_+;\dot B^{n/p-1}_{p,1})} \lesssim 
 \|D\dv\|_{L_2(\R_+;\dot B^{n/p-1}_{p,1})}\quad\hbox{if }\ p<2n.
\end{equation}

 \end{lem}
 \begin{p}
 In order to prove the first inequality, we use the fact
 that, for $i=1,2,$ we have
 $$
 A_i=(\Id+C_i)^{-1}=\sum_{k\geq0}(-1)^kC_i^k\quad\hbox{with}\quad
 C_i(t)=\int_0^t D\bar v_i\,d\tau.
 $$
 Hence
$$
A_2-A_1= \sum_{k\geq 1} \Bigl(C_2^k-C_1^k\Bigr)
=\biggl(\int_0^tD\dv\,d\tau\biggr)\sum_{k\geq1}\sum_{j=0}^{k-1}
C_1^jC_2^{k-1-j}.
$$
So using the fact that $\dot B^{n/p}_{p,1}$ is a Banach algebra, it is easy to conclude
to \eqref{eq:dA}.
\smallbreak
The second inequality is a consequence of Lemma \ref{l:adj} and of the Taylor formula which ensures
that, denoting $\dC:=C_2-C_1,$
$$
\adj(DX_2)-\adj(DX_1)=(\Tr(\dC))\Id-\dC+dP_2(C_1)(\dC)+\frac12d^2P_2(C_1,C_1)(\dC,\dC)+\dotsm
$$
where the coefficients of $P_2$ are polynomials of degree $n-1.$ 
As the sum is finite and $\dot B^{n/p}_{p,1}$ is a Banach algebra, we get \eqref{eq:dAdj}.
\smallbreak
In order to prove the last two estimates, it is only a matter of differentiating the above
relation with respect to $t.$ Keeping in mind the definition of $C_1$ and $C_2,$ we get
$$
\d_t(\adj(DX_2)-\adj(DX_1))=(\div\dv)\Id-D\dv+dP_2(C_1)\cdot\d_t\dC+d^2P_2(C_1)(\d_tC_1,\dC)+\dotsm.
$$
Then using the product laws in Besov spaces yields the desired inequalities.
\end{p}
\medbreak
Finally, we have to justify that the multiplier space $\cM(\dot B^{n/p-1}_{p,1})$
contains characteristic functions of $C^1$ bounded domains, if $p>n-1.$ 
This is a consequence of the following lemma. 
\begin{lem}\label{l:jump}
 Let $\Omega$ be the half-space $\R^n_+$ or a bounded domain  of $\R^n$ with $C^1$ boundary. 
 Assume that  $s\in \R$ and  $p,q \in [1,\infty]$ are such that 
 \begin{equation}\label{eq:trace}
 -1+\frac 1p < s < \frac 1p\cdotp
 \end{equation}
Then the characteristic function 
$\chi_\Omega$ of $\Omega$ belongs to the space ${\mathcal M}(\dot B^s_{p,q}(\R^n)).$
\end{lem}

\begin{p}
This  result which belongs to the mathematical folklore is closely related
to the fact that under Condition \eqref{eq:trace}, functions in $\dot B^s_{p,q}(\Omega)$
extended by $0$ on the whole space, belong to $\dot B^s_{p,q}(\R^n).$
In the case where $\Omega$ is the half-space $\R^n_+$ the
lemma  has been proved in \cite{DM2}, Prop. 3.

If $\Omega$ is a bounded $C^1$ domain, then one may find a finite number $N$
of $\cC^\infty_c(\R^n)$ functions  $\phi_i$ and $C^1$ diffeomorphisms $\psi_i$
so that for any $u\in \dot B^s_{p,q}(\R^n),$ 
$$
u1_\Omega=\sum_{i=1}^N u\phi_i 1_\Omega\quad\hbox{and}\quad
(u\phi_i 1_\Omega)\circ\psi_i= 1_{\R^+}\cdot\bigl((u\phi_i)\circ\psi_i\bigr).
$$
Now, because  Condition \eqref{eq:trace} is satisfied, 
the space $\dot B^s_{p,q}(\R^n)$ is stable by  multiplication 
by smooth compactly supported functions, and by 
$C^1$ change of variables (see e.g. \cite{DM3}, Chap. 2).
Therefore the functions  $(u\phi_i)\circ\psi_i$ belong
to $\dot B^s_{p,q}(\R^n),$ too. 
Using again the stability of this space by multiplication by $1_{\R^+},$ 
one may thus conclude that $(u\phi_i 1_\Omega)\circ\psi_i\in \dot B^s_{p,q}(\R^n).$ Hence
$u\phi_i 1_\Omega$ is in $\dot B^s_{p,q}(\R^n)$ for all $i\in\{1,\cdots,N\}.$ This completes 
the proof of the lemma. 
 \end{p}

\medskip 

\smallskip

{\footnotesize 
\noindent{\bf Acknowledgment.} 
 The second author (PBM) has been partly supported by Polish MN grant No. N N201 547438 and by Foundation for Polish Science in fr. EU European Regional Development Funds (OPIE 2007-2013). He thanks the University Paris-Est Cr\'eteil, where a part of the paper has been performed, for its  kind hospitality. }


\begin{thebibliography}{99}
\bibitem{Abidi} H. Abidi: \'Equation de Navier-Stokes avec densit\'e et viscosit\'e variables dans l'espace critique, {\em Rev. Mat. Iberoam.}, {\bf  23}(2), pages  537--586 (2007). 

\bibitem{AP} H. Abidi and M. Paicu:
Existence globale pour un fluide inhomog\`{e}ne,
  {\it Annales de l'Institut Fourier}, {\bf 57}(3), pages 883--917 (2007).
  
  \bibitem{AKM} S. Antontsev, A. Kazhikhov
and V. Monakhov:  Boundary value problems in mechanics of nonhomogeneous
fluids.  Studies in Mathematics and its
Applications, 22. North-Holland Publishing Co., Amsterdam, 1990.
  
\bibitem{BCD} 
H. Bahouri, J.-Y. Chemin and  R. Danchin: {\it Fourier Analysis and Nonlinear Partial Differential Equations,} Grundlehren der mathematischen Wissenschaften, {\bf 343}, 
Springer (2011).

\bibitem{CK} H.J. Choe and H. Kim: Strong solutions of the Navier-Stokes equations for nonhomogeneous incompressible fluids, {\em Comm. Partial Differential Equations},
{\bf 28} no. 5-6, 1183--1201 (2003).

\bibitem{D1} R. Danchin: Density-dependent incompressible viscous fluids in critical spaces,
{\it Proceedings of the  Royal  Society of  Edinburgh, Sect. A},
{\bf 133}(6), pages 1311--1334 (2003). 

\bibitem{DM2} R. Danchin and P. B. Mucha:
A critical functional framework for the inhomogeneous Navier-Stokes equations in the half-space, {\em J. Funct. Anal.}, {\bf 256}(3), pages 881--927 (2009).

\bibitem{DM3} R. Danchin and P.B. Mucha: Critical functional framework and maximal regularity
in action on systems of incompressible flows, in progress. 

\bibitem{Germain} P. Germain: Strong solutions and weak-strong uniqueness for the nonhomogeneous Navier-Stokes equation, {\em J. Anal. Math.}, {\bf 105}, pages 169--196 (2008).  

\bibitem{Hoff} D. Hoff: Uniqueness of weak solutions of the Navier-Stokes equations of multidimensional
compressible flow, {\it SIAM Journal on Mathematical Analysis}, {\bf 37}(6), pages 1742--1760 (2006). 

\bibitem{LS}
O. Ladyzhenskaya and V. Solonnikov:  The unique solvability of an initial-boundary value problem for viscous incompressible inhomogeneous fluids, {\em  
Journal of Soviet Mathematics}, {\bf  9}, pages 697--749 (1978). 

\bibitem{Lions} P.-L. Lions: Mathematical Topics in Fluid Dynamics,
Vol.~$1$ Incompressible Models, {\it Oxford University Press}  (1996).

\bibitem{MS} 
V. Maz'ya and T. Shaposhnikova:  {\it Theory of Sobolev multipliers. With applications to differential and integral operators.} Grundlehren der Mathematischen Wissenschaften, 
{\bf 337}, Springer (2009).

\bibitem{Mu}
P.B. Mucha:  On weak solutions to the Stefan problem with Gibbs-Thomson correction,
{\em Differential and  Integral Equations}, {\bf  20}(7), 769--792 (2007).
  
\bibitem{MZ}
P.B. Mucha and W.M. Zaj\c{a}czkowski: On local existence of solutions of free boundary problem for incompressible viscous self-gravitating fluid 
motion, {\it Applicationes Mathematicae}, {\bf 27}(3), pages 319--333 (2000). 

\bibitem{So}
V.A. Solonnikov: On the nonstationary motion of isolated value of viscous incompressible fluid,
{\em  Izv. AN SSSR,} {\bf 51}(5), pages  1065--1087 (1987).
\end{thebibliography}
\end{document}